\documentclass{amsart}
\usepackage{amsmath}
\usepackage{amssymb}
\usepackage{graphicx}
\usepackage{amsthm,enumitem}
\newtheorem{theorem}{Theorem}[section]

\newtheorem{lemma}[theorem]{Lemma}

\newtheorem{proposition}[theorem]{Proposition}
\theoremstyle{definition}

\newtheorem{example}{Example}
\newtheorem{remark}{Remark}
\newenvironment{Proof}{{\textit{Proof}.}\ }{~$\square$\vspace{0.2truecm}}

\begin{document}

\title[Locally solvable maximal subgroups]{Locally solvable maximal subgroups\\ in division rings}

\author[Huynh Viet Khanh]{Huynh Viet Khanh}	
\email{khanhhv@hcmue.edu.vn} 
\address{Department of Mathematics and Informatics, HCMC University of Education, 280 An Duong Vuong Str., Dist. 5, Ho Chi Minh City, Vietnam}

\author[Bui Xuan Hai]{Bui Xuan Hai}
\email{bxhai@hcmus.edu.vn} 
\address{Faculty of Mathematics and Computer Science, University of Science, Ho Chi Minh City, Vietnam} 
\address{Vietnam National University, Ho Chi Minh City, Vietnam}


\maketitle

\begin{abstract} 
	Let $D$ be a division ring with center $F$, and $G$ an almost subnormal subgroup of $D^*$. In this paper, we show that if $G$ contains a non-abelian locally solvable maximal subgroup, then $D$ must be a cyclic algebra of prime degree over $F$. Moreover, it is proved that every locally nilpotent maximal subgroup of $G$ is abelian.
\end{abstract}

\section{Introduction and statement of results}
	The existence of a non-cyclic free subgroup in linear groups dates back to 1972 in a ground-breaking article of Jacques Tits \cite{tits}, who proved that every finitely generated subgroup of ${\rm GL}_n(F)$, where $F$ is a field and $n\geq1$, contains either a non-cyclic free subgroup or a solvable normal subgroup of finite index. This result is now commonly referred to as the Tits Alternative for linear groups. As a starting point, Tits' theorem has motivated many efforts in attempts to find its analogous version for skew linear groups; that is, whether Tits' theorem also holds for subgroups of ${\rm GL}_n(D)$, where $D$ is a non-commutative division ring. Unfortunately, the work of Lichtman~\cite{lichtman} has revealed that there is no such an analogy even for the general skew linear group of degree 1; that is, for ${\rm GL}_1(D)$. More precisely, he gave the example of  a division ring whose multiplicative group contains a finitely generated subgroup which possesses neither non-cyclic free subgroups nor solvable normal subgroups of finite index. The question whether the multiplicative group of a non-commutative division ring contains a non-cyclic free subgroup still remains open until now. Meanwhile, it has been shown that the question has a positive answer in some particular cases. For example, a non-commutative division ring contains non-cyclic free subgroups if it is weakly locally finite (\cite{hai-khanh-acta}) or if its center is uncountable (\cite{chiba}). 

	As another direction, some authors also studied the existence problem of non-cyclic free subgroups in maximal subgroups, say $M$, of ${\rm GL}_n(D)$. In the case when $D$ is a centrally finite division ring, it has been proved that Tits Alternative also holds for such maximal subgroups $M$ (see e.g., \cite{mahdavi-free}, \cite{mahdavi-tits} and \cite{moghaddam-mahdavi-17}). Even more, the paper \cite{kiani-mahdavi-2010} demonstrated that this result can be extended to maximal subgroups of a normal subgroup, rather than the whole group ${\rm GL}_n(D)$. Recently, the authors have considered the same problem in a more general situation: the base division ring $D$ is a locally finite dimensional division algebra, and $M$ is assumed to be maximal in an almost subnormal subgroup. Under such assumptions, we have obtained an analogous result in \cite[Theorem 3.1]{hai-khanh}.

	In this paper, we focus our attention on maximal subgroups of the skew linear groups of degree 1. The problem of the existence of maximal subgroups in multiplicative group of a non-commutative division ring is difficult  and unsolved until now. As it was pointed out by  R. Hazrat and A. R. Wadsworth in \cite{ha-wa}, the existence of a non-commutative centrally finite division ring with no maximal subgroups requires in particular the existence of a non-cyclic division algebra of prime degree. The existence of such non-cyclic algebras is exactly Albert's famous conjecture about the cyclic division algebras of prime index.  Over the last decade, group-theoretic  properties of maximal subgroups of the multiplicative group of a division ring have been studied intensively by many authors. The most basic results in this direction can be found in the excellent work of S. Akbari, R. Ebrahimian, H. Momenaee Kermani and A. Salehi Golsefidy (\cite{akbari}). This paper has provided a motivation for many subsequent investigations of various aspects of such maximal subgroups. We refer to the self-contained survey of R. Hazrat, M. Hahdavi-Hezavehi and M. Motiee for details about this area of study (\cite[Section 4]{hazrat}). Let $D$ be a division ring with center $F$, and $D^*=D\backslash\{0\}$ its multiplicative group. It was shown in \cite{ebrahimian_04} that every nilpotent maximal subgroup of $D^*$ is abelian. At the other extreme, if the word `nilpotent' is substituted by `solvable', then the result is no longer true. In fact, in the multiplicative group $\mathbb{H}^*$ of the division ring of real quaternions $\mathbb{H}$, the set $\mathbb{C}^*\cup \mathbb{C}^*j$ is a non-abelian solvable maximal subgroup (\cite{akbari}). More generally, it was proved that if $D^*$ with center $F$ contains a non-abelian solvable maximal subgroup, then $D$ must be a cyclic algebra of prime degree over $F$ (\cite{dorbidi2011}). In \cite{hai-tu}, this result was slightly extended to the case of local solvability. In fact, it was proved that the same result also holds if $D^*$ contains a non-abelian locally solvable maximal subgroup whose derived subgroup is algebraic over $F$.
 
 	Many other authors have generalized the above results to maximal subgroups of a subnormal subgroup $G$, instead of the whole group $D^*$. In particular, M. Ramezan-Nassab and D. Kiani have obtained in \cite{nassab14} an analogous result for nilpotent maximal subgroups of $G$; namely, they have shown that every nilpotent maximal subgroup of $G$ is abelian. Also, the authors and Fallah-Moghaddam independently proved similar result for non-abelian solvable maximal subgroups (see \cite{khanh-hai} and \cite{moghaddam}, respectively). At this point, it is reasonable to ask whether there exist analogous results for locally nilpotent maximal subgroups or locally solvable maximal subgroups of a subnormal subgroup $G$ of $D^*$. In \cite{khanh}, the first author has provided a partial positive answer in the case when $G$ contains a non-abelian locally solvable maximal subgroup, say $M$, such that the $i$-th derived subgroup $M^{(i)}$ of $M$ is algebraic over $F$. In this paper, we shall give a complete affirmative answer to this question without any assumption on $M$. We also extend beyond the above results to the case where $G$ is an almost subnormal subgroup of $D^*$.
 
 	Following Harley \cite{Hartley_89}, a subgroup $A$ of a group $B$ is called \textit{almost subnormal} in $B$ if there is a finite sequence of subgroups
	 $$
	 A=A_0\leq A_1\leq\cdots\leq A_n=B,
	 $$
 	in which either $[A_{i+1}:A_i]$ is finite or $A_i$ is normal in $A_{i+1}$ for $0\leq i\leq n-1$. Clearly if $A$ is a subnormal subgroup of $B$, then it is almost subnormal. But we do not have the converse. In \cite{deo-bien-hai-19} and \cite{ngoc_bien_hai_17}, there are examples of division rings whose multiplicative groups contain almost subnormal subgroups that are not subnormal. 
 
 	The present paper is organized as follows. Section 2 presents examples of maximal subgroups of almost subnormal subgroups. In Section 3, we prove the following main results.

\begin{theorem}\label{theorem_1.1}
	Let $D$ be a division ring with center $F$, and $G$ an almost subnormal subgroup of $D^*$. If $M$ is a locally nilpotent maximal subgroup of $G$,  then $M$ is abelian.
\end{theorem}

\begin{theorem}\label{theorem_1.2}
	Let $D$ be a division ring with center $F$, and $G$ an almost subnormal subgroup of $D^*$. If $M$ is a non-abelian locally solvable maximal subgroup of $G$, then the following hold:
	\begin{enumerate}[font=\normalfont]
		\item[(i)] There exists a maximal subfield $K$ of $D$ such that $K/F$ is a finite Galois extension with $\mathrm{Gal}(K/F)\cong M/K^*\cap G\cong \mathbb{Z}_p$ and $[D:F]=p^2$, for  some prime number $p$. 
		\item[(ii)] The subgroup $K^*\cap G$ is the $FC$-center. Also, $K^*\cap G$ is the Hirsch-Plotkin radical of $M$. For any $x\in M\setminus K$, we have $x^p\in F$ and $D=F[M]=\bigoplus_{i=1}^pKx^i$.
\end{enumerate}
\end{theorem}

	Let $D$ be a non-commutative division ring and $G$ an almost subnormal subgroup of $D^*$. In view of these theorems, one can easily see that every locally solvable maximal subgroup of $G$ is solvable of class at most $2$. Perhaps the best-known example of such a maximal subgroup is the subgroup $\mathbb{C}^*\cup \mathbb{C}j$ of the multiplicative group of the division ring of real quaternion $\mathbb{H}$, and this maximal subgroup is non-abelian. The case of abelian maximal subgroups is very little-known. Until now, there has been no examples of abelian maximal subgroups of the multiplicative group of a non-commutative division ring. It was conjectured in \cite{akbari} that if a division ring contains an abelian maximal subgroup then it is a field. 

	Throughout this paper, for a ring $R$ with the identity $1\ne0$, the symbol $R^*$ stands for the group of units of $R$. Let $S$ be a subring of $R$, and $G$ a subgroup of $R^*$ normalizing $S$ such that $R=S[G]$. Suppose that $N=G\cap S$ is a normal subgroup of $G$ and $R=\bigoplus_{t\in T} tS$, where $T$  is some (and hence any) transversal $T$ of $N$ to $G$. Then, we say that $R$ is a \textit{crossed product} of $S$ by $G/N$ (see \cite{wehrfritz_91} or \cite[p.23]{shirvani-wehrfritz}). 
	
	For a group $G$, let $\tau(G)$ be the unique maximal periodic normal subgroup of $G$, and $B(G)$ a subgroup of $G$ such that $B(G)/\tau(G)$ is the Hirsch-Plotkin radical of $G/\tau(G)$.  For a positive integer $i$, the symbol $G^{(i)}$ denotes the $i$-th derived subgroup of $G$. If $A$ is a ring or a group, then $Z(A)$ denotes the center of $A$. If $D$ is a division ring with the center $F$ and $\Omega$ is a subset of $D$, then $F[\Omega]$ and $F(\Omega)$ denotes respectively the subring and the division subring of $D$ generated by  $F\cup \Omega$.

\section[Examples]{Examples of maximal subgroups}

	The existence problem of maximal subgroups of the multiplicative group of a non-commutative division ring is unsettled completely. An excellent reference for this is the paper of R. Hazrat and A. R. Wadsworth \cite{ha-wa} in which it was shown that ``most of'' centrally finite division rings contain maximal subgroups. In this section, by using the well-known Malcev-Neumann construction, we will provide examples of centrally infinite division rings which contain non-normal maximal subgroups. For this purpose, we first record the following lemma whose proof follows immediately from the correspondence theorem of groups.

\begin{lemma}\label{lemma_2.1}
	Let $G$ and $H$ be groups, and $\varphi: G\to H$ a surjective  group homomorphism. If $N$ is a maximal subgroup of $H$, then $\varphi^{-1}(N)$, the inverse image of $N$, is a maximal subgroup of $G$. Moreover, if $N$ is non-normal in $H$, then $\varphi^{-1}(N)$ is non-normal in $G$.
\end{lemma}

	Now, we shall present here  a brief account of the construction of the general Mal'cev-Newmann series ring, and we refer the reader to  \cite[p. 230]{lam} for more details. Let $R$ be a ring and $G$ an ordered group. Fix a group homomorphism $\omega: G\to {\rm Aut}(R)$ and write $\omega_g$ for the image of $g\in G$ under $\omega$. As a set, the Mal'cev-Newmann ring consists of formal (possibly infinite) sum of form
	$$
	\alpha=\sum_{g\in G}\alpha_gg,
	$$
	where the $\alpha_g$'s are element of $G$. For each such $\alpha$, define the \textit{support} of $\alpha$ by $supp(\alpha):=\{g\in G:\alpha_g\ne0\}$. The Mal'cev-Newmann series ring $A=R((G,\omega))$ is defined to be the set of series over $R$ with well-ordered support; that is,
	$$
	A=R((G,\omega))=\{\alpha=\sum_{g\in G}\alpha_gg: supp(\alpha)\subseteq \text{ is well-orderd}\}.
	$$
	Then, we can make $A=R((G,\omega))$ to be a ring by adding and multiplying its elements according to the following rules.
	$$
	\sum_{g\in G}\alpha_gg+\sum_{g\in G}\beta_gg=\sum_{g\in G}(\alpha_g+\beta_g)g,
	$$
	
	$$
	\left( \sum_{g\in G}\alpha_gg\right) \left( \sum_{h\in G}\beta_hh\right) =\sum_{t\in G}\left( \sum_{gh=t}\alpha_g\omega_g(\beta_h)\right) t.
	$$
	The multiplication of two series in $A=R((G,\omega))$ is built from the twist law 
	\[
	g\cdot r=\omega_g(r)g \tag{$\star$}
	\]
	where $r\in R$ and $g\in G$. It is shown in \cite[Theorem 14.21]{lam} that $A=R((G,\omega))$ is a division ring whenever $R$ is a division ring.
	
	With suitable choices of $R$, $G$ and $\omega_g$, one can obtain from $A=R((G,\omega))$ centrally finite division rings. For example, when $k:=R$ is a field, $G$ is an infinite cyclic group $\{x^n: n\in\mathbb{Z}\}$ and $\omega$ is specified by a single automorphism $\sigma:=\omega_x$, the twist law ($\star$) boils  down to $x\cdot r=\sigma(r)x$ for all $r\in R$. In this case, the Mal'cev-Newmann series ring $k((\left\langle x\right\rangle ,\omega))$ is exactly the Hilbert's twisted Laurent division ring $k((x,\sigma))$. And, it is shown in \cite[Proposition 14.2]{lam} that $k((x,\sigma))$ is centrally finite if and only if $\sigma$ has a finite order. As another example, in \cite[Section 2]{tignol-amitsur}, using the Mal'cev-Newmann construction with a suitable co-cycle $f$, J.-P. Tignol and S. A. Amitsur has constructed the division ring $D_f$ which is finite dimensional over its center $C_f$ equipped with a valuation $v: D_f^*\to\mathbb{Z}^n$ with value group $\mathbb{Z}^n$. 
	
	There are also some ways to obtain centrally infinite division rings from Mal'cev-Newmann series ring $A=R((G,\omega))$. For instance, when $R$ is a field, $G$ is a non-trivial ordered group and $\omega:G\to{\rm Aut}(R)$ is an injective homomorphism, then $Z(A)=\{r\in R: \omega_g(r)=r \text{ for all } g\in G\}$ and the division ring $A$ is centrally infinite (see \cite[Corollary 14.26]{lam}). 
	
	Now, we give another way to obtain centrally infinite division ring in which we will construct maximal subgroups of almost subnormal subgroups using Mal'cev-Newmann construction $R((G,\omega))$. For this purpose, take $K:=R$ to be a field, $G$ to be a free group generated by $\{x_i:i\in I\}$ with $|I|\geq2$ and $\omega_g$ to be the trivial homomorphism. Since any free group can be ordered (\cite[Theorem 6.31]{lam}), the Mal'cev-Newmann series ring $D:=K((G,\omega))$ is a division ring.  As $\omega$ is trivial, the twist law ($\star$) above can be written as $g \cdot r = rg$ for all $r\in R$ and $g\in G$; that is, the elements of $G$ commute elementwise with the element of $R$. Also, we may write $D=K((G)$ for $D:=K((G,\omega))$ due to the triviality of $\omega$. It is shown in \cite[Corollary 14.25]{lam} that the free ring $K\left\langle x_i:i\in I\right\rangle $ can be embedded in $D=K((G))$. Since $G$ is free of rank $\geq 2$, one can verify that the center of $D$ is  $K$, and $D$ is infinite dimensional over $K$. Moreover, the function $v:D^*\to G$ defined by $\alpha\mapsto\min\{supp(\alpha)\}$ is a surjective group homomorphism. (This is a special case of Exercise 14.10, p.236 of \cite{lam} whose solution can be found in \cite{lam-ex}.)

\begin{example}
	Let $H$ be a finite group of order $pq$, where $p$ and $q$ are distinct prime numbers. By Sylow's Theorem, $H$ contains a subgroup $N$ of order $p$ (and thus of index $q$). Moreover, there exist a free  subgroup $G$ of rank $\geq2$ and a surjective group homomorphism $u: G\to H$. Thus, for an arbitrary field $K$ we can construct the Malcev-Neumann (centrally infinite) division ring $D=K((G))$ as well as a surjective $v:D^*\to G$ as above. If we set $\varphi=u\circ v$, then $\varphi:D^*\to H$ is a surjective homomorphism and $D^*/\ker\varphi\cong H$. Setting $A=\varphi^{-1}(N)$  and $B=\ker\varphi$, we obtain that $B\leq A\leq D^*$ with $[D^*:A]=q$ and $[A:B]=p$. It is clear that $A$ is an almost subnormal maximal subgroup of $D^*$, and $B$ is a maximal subgroup of $A$. Moreover, if $N$ is chosen to be non-normal in $H$, then it follows from Lemma \ref{lemma_2.1} that $A$ is not normal in $D^*$.
\end{example}

	We note in passing that, if $H$ is chosen to be infinite rather than finite, then one can apply Lemma \ref{lemma_2.1} to produce more interesting maximal subgroups of $D^*$. 

\section[Locally solvable maximal subgroups]{Locally solvable maximal subgroups}
	This section is devoted to examining locally solvable maximal subgroups of an almost subnormal subgroup in a division ring.  For this purpose, we prepare the way by first establishing a few group-theoretic lemmas which are easily proved.

\begin{lemma}\label{lemma_3.1} 
	Let $H\leq G$ be groups. If $H$ is a characteristic subgroup of $G$, then so is $C_G(H)$. 
\end{lemma}

\begin{lemma}\label{lemma_3.2} 
	Let $A\leq H\leq G$ be groups. If $A$ and $H/A$ are characteristic subgroups of $G$ and $G/A$, respectively, then $H$ is  a characteristic subgroup of $G$.
\end{lemma}

\begin{lemma}\label{lemma_3.3}
	Every locally solvable periodic group is locally finite.
\end{lemma}
\begin{Proof}
	To begin with, we assume that $G$ is a locally solvable periodic group. If $H$ is any finitely generated subgroup of $G$, then our hypothesis implies that $H$ is solvable of derived length $n$, say. In other words, we have a series
	$$1=H^{(n)}\unlhd H^{(n-1)}\unlhd\cdots\unlhd H'\unlhd H.$$
	We shall prove that $H$ is finite by induction on its derived length $n$. If $n=1$, then our conclusion follows from the fact that every finitely generated periodic abelian group is finite. We may therefore assume that $n>1$. Since $H/H'$ is abelian as well as periodic finitely generated, it is also finite. In other words, $H'$ is a subgroup of finite index in the finitely generated group $H$. This fact assures us that $H'$ is also finitely generated. By induction hypothesis, we conclude that $H'$ is finite and, in consequence, $H$ is finite too.
\end{Proof}

\begin{lemma}\label{lemma_3.4}
	Let $D$ be a division ring with center $F$, and $G$ an almost subnormal subgroup of $D^*$. If $G$ is locally solvable, then it is contained in $F$.
\end{lemma}
\begin{Proof}
	Assume by contradiction that $G$ is not contained in $F$. It follows from \cite[Proposition 2.2]{khanh-hai-22} that $G$ contains a non-central subnormal subgroup which is also locally solvable. This contradicts the fact that every locally solvable subnormal subgroup of $D^*$ is central (\cite[Theorem 2.6]{2019_danh-khanh}). The proof is therefore proved. 
\end{Proof}

	The next two lemmas could have been proved by some other authors, for completeness we include their proofs here, since we were unable to find ones in the literature.
\begin{lemma}\label{lemma_3.5} 
	Let $D$ be a division ring with center $F$, and $G$ a subgroup of $D^*$. If $G$ is locally finite, then $F(G)=F[G]$.  
\end{lemma}
\begin{Proof}
	Pick any non-zero element $x$ in $F[G]$, we need only prove that $x^{-1}$ also belongs to $F[G]$. By the choice of $x$, we can find some elements $f_1,f_2,\dots, f_n$ in $F$ and $g_1,g_2,\dots, g_n$ in $G$ such that
	$$x=f_1g_1+f_2g_2+\cdots +f_ng_n.$$
	Let $H$ be the subgroup of $G$ generated by the $g_i$'s. Then, the local finiteness of $G$ ensures that $H$ is a finite group. This obviously implies that $F[H]$ is a domain which is a finite dimensional vector space over $F$. Accordingly, this is then actually a division ring and so $F(H)=F[H]$. Consequently, we have $x^{-1}\in F[H]\subseteq F[G]$. 
\end{Proof}

\begin{lemma}\label{lemma_3.6}
	Let $D$ be a division ring with center $F$, and $G$ a locally finite subgroup of $D^*$. If $F(G)=D$, then $D$ is a locally finite dimensional division algebra.
\end{lemma}
\begin{Proof}
	It follows from Lemma \ref{lemma_3.5} that $D=F[G]$. Let $\{x_1,x_2,\dots,x_k\}$ be a finite subset of $ D$. Then, for each $i$, there exist elements $f_{i_1}, f_{i_2},\dots, f_{i_s}$ in $F$ and $g_{i_1}, g_{i_2},\dots, g_{i_s}$ such that
	$$x_i=f_{i_1}g_{i_1}+f_{i_2}g_{i_2}+\cdots+f_{i_s}g_{i_s}.$$
	Set $H$ to be the subgroup of $G$ generated by all the $g_{i_j}$'s.
	We know that $G/G\cap F^*\cong GF^*/F^*$, which is clearly a locally finite group. This facts yields that $HF^*/F^*$ is finite; thus, we can take a transversal $\{y_1,y_2,\dots,y_t\}$ of $F^*$ in $HF^*$. Moreover, if we set
	$$S=Fy_1+Fy_2+\cdots+Fy_t,$$ 
	then this is a domain containing $\{x_1,x_2,\cdots,x_k\}$. Besides, the above relation reveals that $S$ is finite dimensional over $F$, so it is actually a division ring. Now, the division subring $F(x_1,x_2,\cdots,x_k)$ is contained in $S$, whence $[F(x_1,x_2,\cdots,x_k):F]$ is a finite. Since the set $\{x_1,x_2,\dots,x_k\}$ is taken to be arbitrary, we conclude that $D$ is locally finite. Our proof is now completed.
\end{Proof}

\begin{lemma}\label{lemma_3.7}
	Let $D$ be a division ring with center $F$, and $G$ an almost subnormal subgroup of $D^*$. Assume that $M$ is a non-abelian locally solvable maximal subgroup of $G$. If $A \unlhd M$, then either $A$ is abelian or $F(A)=D$.
\end{lemma}
\begin{Proof}
	Since $A \unlhd M$, we have $M\leq G\cap N_{D^*}(F(A)^*)\leq G$. The maximality of $M$ in $G$ implies that either $G\cap N_{D^*}(F(A)^*)=M$ or $G\leq N_{D^*}(F(A)^*)$. If the first case occurs, then $A \unlhd F(A)^*\cap G$, so $A$ is almost subnormal in $F(A)^*$ contained in $M$. Since $M$ is locally solvable, so is $A$. It follows from Lemma \ref{lemma_3.4} that $A$ is contained in the center of $F(A)$, which implies that $A$ is abelian. In the latter case, the division subring $F(A)$ of $D$  is normalized by $G$. With reference to \cite[Proposition 2.2]{khanh-hai-22}, we conclude that $G$ contains a non-central subnormal subgroup that also normalizes $F(A)$. It follows from \cite[Theorem 1]{stuth} that either $A\subseteq F$ or $F(A)=D$. Hence, either $A$ is abelian or $F(A)=D$. Our proof is now finished.
\end{Proof}

For the convenience of readers, we gather a number of theorems of B. A. F Wehrfritz, which will be used in the proofs of our results.

\begin{lemma}[{\cite[2.5]{wehrfritz86_loc_nil_I}}]\label{lemma_3.8}
	Let $R=F[G]$ be an $F$-algebra, where $F$ is a field and $G$ is a locally nilpotent group of units of $R$, such that for every finite subset $X$ of $R$ there is a finitely generated subgroup $Y$ of $G$ with $F[Y]$ prime and containing $X$. Let $\tau(G)$ be the unique maximal periodic normal subgroup of $G$, and $Z/\tau(G)$ the center of $G/\tau(G)$. Then, $R$ is a crossed product of $F[Z]$ by $G/Z$.
\end{lemma}

\begin{lemma}[{\cite[7]{wehrfritz88_local_nil_II}}]\label{lemma_3.9}
	Let $R=F[G]$ be an $F$-algebra, where $F$ is a field and $G$ is a locally solvable subgroup of the group of units of $R$ such that for every infinite subgroup $X$ of $G$ the left annihilator of $X-1$ in $R$ is $\{0\}$. Let $B(G)$ be a subgroup of $G$ such that $B(G)/\tau(G)$ is the Hirsch-Plotkin radical of $G/\tau(G)$. Then $R$ is a crossed product of $F[B(G)]$ by $G/B(G)$.
\end{lemma}

\begin{lemma}[{\cite[3.2]{wehrfritz_91}}]\label{lemma_3.10}
	Let $R$ be a ring, $J$ a subring of $R$, and $H\leq K$ subgroups of the group of units of $R$ normalizing $J$ such that $R$ is the ring of right quotients of $J[H]\leq R$ and $J[K]$ is a crossed product of $J[B]$ by $K/B$ for some normal subgroup $B$ of $K$. Then $K=HB$.
\end{lemma}

For a group $G$, the second center $Z_2(G)$ of $G$ is defined by $$Z_2(G)/Z(G)=Z(G/Z(G)).$$

\begin{lemma}[{\cite[Theorem 1.1(c)]{wehrfritz_91}}]\label{lemma_3.11}
	Let $G$ be a $\left\langle P, L\right\rangle\mathfrak{A} $-subgroup of ${\rm GL}_n(D)$ such that the subalgebra $F[N]$ of ${\rm M}_n(D)$ is a prime ring for every characteristic subgroup $N$ of $G$. Denote by $\tau(G)$ the unique maximal locally finite normal subgroup of $G$. If $\tau(G)\subseteq Z_2(G)$, then $F[G]$ is a crossed product of $F[A]$ by $G/A$, for some abelian characteristic subgroup of $G$.
\end{lemma}

\begin{remark}\label{remark}
	In the above lemma, the notation $\left\langle P, L\right\rangle\mathfrak{A} $ stands for a class of groups, in which $\mathfrak{A} $ denotes the class of abelian groups, $P$ and $L$ the poly and local operators (\cite{wehrfritz_91}). It is a simple matter to check that the class of groups $\left\langle P, L\right\rangle\mathfrak{A} $ contains that of locally solvable groups.
\end{remark}

The following proposition is a simple modification of  \cite[Theorem 3.3]{hai-tu}, so it should be omitted.

\begin{proposition}\label{proposition_3.12}
	Let $D$ be a division ring with center $F$, and $G$ an almost subnormal subgroup of $D^*$. If $M$ is a non-abelian metabelian maximal subgroup of $G$, then $[D:F]<\infty$.  
\end{proposition}

\begin{lemma}[{\cite[Proposition 4.1]{wehrfritz_07}}]\label{lemma_3.13}
	Let $D=E(A)$ be a division ring generated by its metabelian subgroup $A$ and its division subring $E$ such that $E\leq C_D(A)$. Set $H=N_{D^*}(A)$, $B=C_A(A')$, $K=E(Z(B))$, $H_1=N_{K^*}(A)=H\cap K^*$, and let $\tau(B)$ be the unique maximal periodic normal subgroup of $B$.
	\begin{enumerate}[font=\normalfont]
		\item[(i)] If $\tau(B)$ has a quaternion subgroup $Q=\left\langle i,j\right\rangle $ of order $8$ with $A=QC_A(Q)$, then $H=Q^+AH_1$, where $Q^+=\left\langle Q,1+j,-(1+i+j+ij)/2\right\rangle$. Also, $Q$ is normal in $Q^+$ and $Q^+/{\left\langle -1,2\right\rangle}\cong Aut Q\cong Sym(4)$.	
		\item[(ii)] If $\tau(B)$ is abelian and contains an element $x$ of order $4$ not in the center of $B$, then $H=\left\langle x+1\right\rangle AH_1$.
		\item[(iii)] In all other cases, $H=AH_1$.
	\end{enumerate}
\end{lemma}

\begin{theorem}\label{theorem_3.14}
	Let $D$ be a division ring with center $F$, and $G$ an almost subnormal subgroup of $D^*$. If $M$ is a non-abelian solvable maximal subgroup of $G$, then the following hold:
	\begin{enumerate}[font=\normalfont]
		\item[(i)] There exists a maximal subfield $K$ of $D$ such that $K/F$ is a finite Galois extension with $\mathrm{Gal}(K/F)\cong M/K^*\cap G\cong \mathbb{Z}_p$ and $[D:F]=p^2$,  for  some prime number $p$. 
		\item[(ii)] The subgroup $K^*\cap G$ is the $FC$-center. Also, $K^*\cap G$ is the Hirsch-Plotkin radical of $M$. For any $x\in M\setminus K$, we have $x^p\in F$ and $D=F[M]=\bigoplus_{i=1}^pKx^i$.
\end{enumerate} 
\end{theorem}
\begin{Proof}
	First, we prove that $[D:F]<\infty$. Since $M$ is non-abelian, Lemma \ref{lemma_3.7} says that $F(M)=D$. Also, we may suppose that $M$ is solvable with derived length $s\geq2$. In other words, we have such a series
	$$1=M^{(s)}\unlhd M^{(s-1)}\unlhd M^{(s-2)}\unlhd\cdots\unlhd M'\unlhd  M.$$ 
	If we set $A=M^{(s-2)}$, then $A$ is a non-abelian metabelian normal subgroup of $M$. As $A$ is non-abelian, Lemma \ref{lemma_3.7} again says that $F(A)=D$, from which it follows that $Z(A)=F^*\cap A$ and $F=C_D(A)$. Set $H=N_{D^*}(A)$,  $B=C_A(A')$, $K=F(Z(B))$, $H_1=H\cap K^*$, and $\tau(B)$ to be the maximal periodic normal subgroup of $B$.
	Then, it is a simple matter to check that $H_1$ is an abelian group. Since $\tau(B)$ is a characteristic subgroup of $B$ (\cite[Lemma 2.1]{2019_danh-khanh}), it follows from Lemma \ref{lemma_3.1} that $\tau(B)$ is characteristic in $A$. In order to use Lemma \ref{lemma_3.13}, we divide our situation into three cases:
	
	\bigskip
	
	\textit{Case 1:} $\tau(B)$ is not abelian. 
	
	\bigskip
	
	Since $\tau(B)$ is characteristic in $B$ and $B$ is normal in $M$, we conclude that $\tau(B)$ is normal in $M$. By virtue of Lemma \ref{lemma_3.7}, we have  $F(\tau(B))=D$. In addition, as $\tau(B)$ is solvable and periodic, it is actually a locally finite group (Lemma \ref{lemma_3.3}). It follows from Lemma \ref{lemma_3.6} that $D=F(\tau(B))=F[\tau(B)] $ is a locally finite dimensional division algebra. Since  $M$ is solvable, it contains no non-cyclic free subgroups. With reference to \cite[Theorem 3.1]{hai-khanh}, we deduce that $[D:F]<\infty$.
	
	\bigskip	
	
	\textit{Case 2:} $\tau(B)$ is abelian and contains an element $x$ of order $4$ not in the center of $B=C_A(A')$.
	
	\bigskip
	
	Since  $x\not\in Z(B)$, it does not belong to $F$. Since $x$ is of finite order, the field $F(x)$ is an algebraic extension of $F$.  Note that $\left\langle x\right\rangle$ is indeed a $2$-primary component of $\tau(B)$ (see \cite[Theorem 1.1, p.132]{wehrfritz_07}); thus, it is a characteristic subgroup of $\tau(B)$. Consequently, $\left\langle x\right\rangle$ is a normal subgroup of $M$. This being the case, all elements of the set $x^M=\{m^{-1}xm\vert m\in M\}\subseteq F(x)$ have the same minimal polynomial over $F$. As a result, $x$ is an FC-element and so $[M:C_M(x)]<\infty$. Now, if we set $C={\rm core}_M(C_M(x))$, then $C$ is a normal subgroup of finite index in $ M$. In view of Lemma \ref{lemma_3.7}, either $F(C)=D$ or $C$ is abelian. The first case cannot occur since it would imply that $x\in F$, which is impossible. Therefore $C$ is abelian. If we set $K=F(C)$, then the finiteness of $M/C$ implies that $K$ is a subfield of $D$, over which $D$ is finite dimensional. This fact yields that $[D:F]<\infty$.
	
	\bigskip	
	
	\textit{Case 3:} $H=AH_1$.
	
	\bigskip
	
	The fact $A'\leq H_1\cap A$ implies that $H/H_1\cong A/A\cap H_1$ is abelian and so $H'\leq H_1$. Since $H_1$ is abelian, it follows that $H'$ is abelian too. Because $M\leq H$, we know that $M'$ is also abelian. In other words, $M$ is a metabelian group; hence, the conclusions follow from Proposition \ref{proposition_3.12}.
	
	By what we have proved, we conclude that $n:=[D:F]<\infty$. Since $M$ is solvable, it contains no non-cyclic free subgroups. In view of \cite[Theorem 3.1]{hai-khanh}, we have $F[M]=D$, there exists a maximal subfield $K$ of $D$  containing $F$ such that $K/F$ is a Galois extension, $N_G(K^*)=M $, $K^*\cap G$ is the Fitting subgroup of $M$ and it is the $FC$-center, and $M/K^*\cap G\cong{\rm Gal}(K/F)$ is a finite simple group of order $[K:F]$.

	Since $M/K^*\cap G$ is solvable and simple, one has $M/K^*\cap G\cong{\rm Gal}K/F)\cong \mathbb{Z}_p$ for some prime number $p$ (\cite[12.5.2, p.367]{robinson}). Therefore, $[K:F]=p$ and $[D:F]=p^2$. For any $x\in M\backslash K$, if $x^p\not\in F$, then by the fact $F[M]=D$, we conclude that $C_M(x^p)\ne M$. Moreover, since $x^p\in K^*\cap G$, it follows that $\left\langle x,K^*\cap G\right\rangle \leq C_M(x^p)$. In other words, $C_M(x^p)$ is a subgroup of $M$ strictly containing $K^*\cap G$. Because $M/K^*\cap G$ is simple, we have $C_M(x^p)= M$, a contradiction. Therefore $x^p\in F$. Furthermore, since $x^p\in K$ and $[D:K]_r=p$, we conclude $D=\bigoplus_{i=1}^{p-1}Kx_i$.  
	
	It remains to prove that $K^*\cap G$ is the Hirsch-Plotkin radical of $M$. Note that we are in the case that $K^*\cap G \subseteq M$, we conclude that $K^*\cap G = K^*\cap M$. Thus, we have $M/K^*\cap M \cong MK^*/K^*\cong \mathbb{Z}_p$. Let $\mathcal{H}$ be the Hirsch-Plotkin radical of $M$. It is clear that $\mathcal{H}K^*/K^*\leq MK^*/K^*$, which is a group of order $p$. Therefore, either $\mathcal{H}\leq K^*$ or $\mathcal{H}K^*=MK^*$. The first case implies that $\mathcal{H}\leq K^*\cap G$, so we actually have $\mathcal{H}= K^*\cap G$; we are done. If the second case occurs, then  $F[\mathcal{H}K^*]=F[M]=D$.  It follows that $\mathcal{H}K^*$ is a locally nilpotent absolutely irreducible subgroup of $D^*$, hence it is center-by-locally finite by \cite[Theorem 1]{wehrfritz88_local_nil_II}. It is a simple matter to check that the center of $\mathcal{H}K^*$ is contained in $F^*$. As a result, $\mathcal{H}K^*/F^*$ is locally finite, from which it follows that $K/F$ is a non-trivial radical Galois extension. According to \cite[15.13]{lam}, we conclude that $D$ is algebraic over a finite subfield. But then Jacobson's Theorem (\cite[13.11]{lam}) says that $D$ is commutative, a contradiction. 
\end{Proof}

\bigskip

	\noindent{\it \textbf{Proof of Theorem~\ref{theorem_1.1}}}. Assume by contradiction that $M$ is non-abelian. The crux of our argument is to show that $n:=[D:F]<\infty$, for then we know that $D\otimes_F D^{op}\cong {\rm M}_n(D)$. Accordingly, we may view $M$ as a subgroup of ${\rm GL}_n(F)$ to conclude that it is actually a solvable group. Then, the previous theorem says that $M$ is distinguish from its Hirsch-Plotkin radical, which is the desired contradiction since $M$ is assumed to be locally nilpotent. 

	Now, we start our argument by setting $\tau(M)$ to be the unique maximal periodic (hence locally finite) normal subgroup of $M$, and $Z/\tau(M)$ to be the center of $M/\tau(M)$. Since $\tau(M)$ is normal in $M$, it follows from Lemma \ref{lemma_3.7} that either $F(\tau(M))=D$ or $\tau(M)$ is abelian. If the first case occurs, then Lemma \ref{lemma_3.6} says that $D$ is a locally finite dimensional division algebra. According to \cite[Theorem~ 3.1]{hai-khanh}, we have $[D:F]<\infty$; we are done. Now assume that $\tau(M)$ is abelian. This fact implies that $Z$ is a solvable group with the derived series length of, say $s\geq 1$. By Lemma ~\ref{lemma_3.2}, we obtain that  $Z$ is a characteristic subgroup of $M$. Again, Lemma~\ref{lemma_3.7} implies  that either $F(Z)=D$ or else $Z$ is abelian. Let us consider the following two possible cases:

	\bigskip
	
	\textit{Case 1:  $F(Z)=D$}. 
	
	\bigskip
	
	Because $D$ is non-commutative, it follows that $Z$ is non-abelian. Therefore, we may suppose that it is solvable with derived length $s\geq2$. As a result, there exists such a series 
	$$1=Z^{(s)}\unlhd Z^{(s-1)}\unlhd Z^{(s-2)}\unlhd\cdots\unlhd Z'\unlhd  Z \unlhd M.$$
	If we set $A=Z^{(s-2)}$, then this is a non-abelian metabelian normal subgroup of $M$. By virtue of Lemma \ref{lemma_3.7}, we must have  $F(A)=D$. It follows that $Z(A)=F^*\cap A$ and $F=C_D(A)$. Set $H=N_{D^*}(A)$,  $B=C_A(A')$, $K=F(Z(B))$, $H_1=H\cap K^*$, and $\tau(B)$ to be the unique maximal periodic normal subgroup of $B$. Then $H_1$ is an abelian group, and $\tau(B)$ is a characteristic subgroup of $B$, and hence of $A$ (Lemma ~\ref{lemma_3.1}). Note that since $M$ is locally nilpotent, it is well-known that $\tau(M)$ is the set of all periodic elements of $M$ and so $\tau(B)\subseteq \tau(M)$, from which it follows that $\tau(B)$ is abelian. Also, it is clear that $\tau(B)$ is a normal subgroup of $M$. Therefore, the same arguments used in Case 2 and Case 3 of the proof of Theorem \ref{theorem_3.14} show that $[D:F]<\infty$.
	
	\bigskip 
	
	\textit{Case 2: }$Z$ is abelian. 
	
	\bigskip 
	
	Let $N$ be the maximal subgroup of $M$ with respect to the property: $N$ is an abelian normal subgroup of $M$ containing $Z$. We assert that $M/N$ is a simple group. For, let $P$ be a normal subgroup of $M$ properly containing $N$. The maximality of $N$ allows us to assume that $P$ is non-abelian, and so $F(P)=D$ by Lemma \ref{lemma_3.7}. Since $P$ is locally nilpotent, it follows from \cite[Corollary 24]{wehrfritz_89_GoldieSubring} that $F[P]$ is an Ore domain whose skew field of fractions is obviously coincided with $D$. With reference to Lemma \ref{lemma_3.8}, we deduce that $F[M]$ is a crossed product of $F[Z]$ by $M/Z$. This being the case, Lemma \ref{lemma_3.10} immediately implies that $M=PZ=P$; recall that $Z\leq P$. This fact shows that $M/Z$ is simple, as asserted. Note that $M\ne Z$ since $M$ is non-abelian. Now, $M/Z$ is both locally nilpotent and simple, we know that it must be finite of prime order. Now, if we set $L=F(Z)$, then $L$ is a subfield of $D$, over which $D$ is finite dimensional. A consequence of this fact is that $[D:F]<\infty$. 
	
	At any rate, we always have $[D:F]<\infty$, which completes our argument. $\square$

	\bigskip

	\noindent{\it \textbf{Proof of Theorem~\ref{theorem_1.2}}}. It suffices to prove that $[D:F]<\infty$, for then, as we have just seen, the group $M$ may be viewed as a linear group. It follows that $M$ is solvable and all results follow from Theorem \ref{theorem_3.14}.  

	Since $M$ is non-abelian, Lemma \ref{lemma_3.7} says that $F(M)=D$. Let $\tau(M)$ be the unique maximal periodic  normal subgroup of $M$. If $\tau(M)$ is non-abelian, then the same arguments used in the beginning of the proof of Theorem \ref{theorem_1.1} shows that $[D:F]<\infty$, finishing our proof. Otherwise, $\tau(M)$ is abelian, and we have two possible cases.
	
	\bigskip 
	
	\textit{Case 1: }$\tau(M)\not\subseteq F$.
	
	\bigskip 

	From the field theory, the field $F(\tau(M))$ is an algebraic extension of $F$. Take $x\in \tau(M)\backslash F$, and set $x^M:=\{m^{-1}xm\vert m\in M\}$. Since $F(\tau(M))$ is normalized by $M$, we have $x^M\subseteq  F(\tau(M))$. Thus, all elements of the set $x^M $ have the same minimal polynomial over $F\subseteq F(\tau(M))$; hence, $x^M$ is finite. In other words, $x$ is an $FC$-element of $M$. By the same argument used in Case 2 of the proof of Theorem ~\ref{theorem_3.14}, we conclude that $[D:F]<\infty$.
	
	\bigskip 

	\textit{Case 2: }$\tau(M)\subseteq F$. 

	\bigskip 
	
	Let $B(M)$ be the subgroup of $M$ such that $B(M)/\tau(M)$ is the Hirsch-Plotkin radical of the group $M/\tau(M)$. Our next step is to  assert that $B(M)$ is indeed a locally nilpotent group. For this purpose, we take an arbitrary finitely generated subgroup $H$ of $B(M)$, and our aim is to show that this is a nilpotent group. It is a simple matter to see that $H\tau(M)/\tau(M)$ is a finitely generated subgroup of $B(M)/\tau(M)$. Accordingly, the local nilpotence of $B(M)/\tau(M)$ implies that $H\tau(M)/\tau(M)$ is nilpotent. We set
	$$H_1=[H,H],\;\;\; H_2=[H_1, H] ,$$
	$$H_3=[H_{2}, H], \;\;\; \cdots\;\;\;\;\;\;\;\;\;\;\;\;\;\;\;\;$$ 
	where $[H,K]$, in particular, stands for the subgroup of $G$ generated by the set of commutators $\{[a,b]=a^{-1}b^{-1}ab|\mbox{ for all }a\in H \mbox{ and } b\in K\}$.
	Now, as $H\tau(M)/\tau(M)$ is nilpotent, we can find an integer $n$ for which $H_n\leq \tau(M)\subseteq F$. This fact says that any element of $H_n$ commutes elementwise with $H$ and, in consequence, we have $H_{n+1}=[H_n,H]=1$, from which it follows that $H$ is nilpotent. In other words, we obtain that $B(M)$ is locally nilpotent, as asserted. It follows that $B(M)$ is the Hirsch-Plotkin radical of $M$. If $M=B(M)$, then $M$ is locally nilpotent and it is thus abelian by Theorem \ref{theorem_1.1}, a contradiction. We may therefore assume that $M\ne B(M)$. 
	
	We claim that $M/B(M)$ is a simple group. For, let $C$ be a normal subgroup of $M$ properly containing $B(M)$. 
	It is obvious that $C$ is non-abelian and so that $F(C)=D$ by  Lemma \ref{lemma_3.7}. With reference to  Lemma \ref{lemma_3.9}, we conclude that $F[M]$ is a crossed product of $F[B(M)]$ by $M/B(M)$. It follows from \cite[Corollary 24]{wehrfritz_89_GoldieSubring} that $F[C]$ is an Ore domain whose skew field of fractions is coincided with $D$. Accordingly, we may apply Lemma \ref{lemma_3.10} to conclude that $M=B(M)C$, which implies that $M=C$. The final fact says that $M/B(M)$ is  simple, as claimed.  
	
	Now, as  $M/B(M)$ is both locally solvable and simple, we obtain that $M/B(M)$ is a finite group of prime order.  Again by Lemma \ref{lemma_3.7}, either $B(M)$ is abelian or else $F(B(M))=D$. Let us consider the following two subcases.
	
	\bigskip
	
	\textit{Subcase 2.1:} $B(M)$ is abelian.
	
	\bigskip
	
	If we set $L=F(B(M))$, then $L$ is a subfield of $D$. Since $F(M)=D$ and $|M/B(M)|<\infty$, we have $[D:L]_r<\infty$. This implies that $[D:F]<\infty$, and we are done.
	
	\bigskip
	
	\textit{Subcase 2.2:} $F(B(M))=D$.
	
	\bigskip
	
	Since we are in the case $\tau(M)\leq F^*\cap M\leq Z_2(M)$, we may use Lemma \ref{lemma_3.11} as well as Remark \ref{remark} to deduce that $F[M]$ is a crossed product of $F[A]$ by $M/A$ for some abelian characteristic subgroup $A$ of $M$. Since $A\tau(M)/\tau(M)$ is an abelian normal subgroup of $M/\tau(M)$, we have $A\tau(M)/\tau(M)\subseteq B(M)/\tau(M)$; recall that $B(M)/\tau(M)$ is the Hirsch-Plotkin radical of $M/\tau(M)$. As a result, we have $A\subseteq B(M)$. By Lemma \ref{lemma_3.10}, we conclude that $M=B(M)A=B(M)$,  contrasting to the fact that $B(M)$ is a proper subgroup of $M$. Finally, we always have $[D:F]<\infty$, which completes our argument.
	
	Before closing this paper, we note that in \cite{khanh-hai-22} the authors also studied locally solvable maximal subgroups of an almost subnormal subgroup of ${\rm GL}_n(D)$, where $n\geq 2$. This case required a different approach by which it was proved that such maximal subgroups are always abelian. 
	
	\textbf{Acknowledgements.}
	The authors are profoundly grateful to the referee for his/her careful review and suggestions which led a significantly improvement to our manuscript.

\end{document}